\documentclass[12pt]{article}
\usepackage{epsfig}
\usepackage{amsmath}
\usepackage{amssymb}
\newtheorem{theorem}{Theorem}

\newtheorem{remark}{Proposition}
\newtheorem{corollary}{Corollary}
\begin{document}

\title{On location of discrete spectrum for complex Jacobi matrices}
\author{I. Egorova, L. Golinski}
\maketitle
\begin{abstract}
We study spectrum inclusion regions for complex Jacobi
matrices that are compact perturbations of the discrete
laplacian. The condition sufficient for the lack of
discrete spectrum for such matrices is given.
\end{abstract}

 {\bf Introduction}.
Let

\begin{equation}
    \label{1.3}
    J=\left(\begin{array}{ccccc}
                b_1&c_1&&&\\
                a_1&b_2&c_2&&\\
                &a_2&b_3&c_3&\\
                &&\ddots&\ddots&\ddots\\
            \end{array}\right)
\end{equation}
be an infinite Jacobi matrix with complex entries. We assume that
$a_nc_n\neq 0$, $n\in\mathbb{N}:=\{1,2,...\},$ and
$$ \lim_{n\to \infty}a_n=\lim_{n\to \infty}c_n=1,\qquad
\lim_{n\to \infty}b_n = 0, $$
that is, the operator $J$ generated by matrix (\ref{1.3})
in $\ell^2( \mathbb{N})$, is a compact perturbation of the
discrete laplacian $$ J_0:\quad a_n=c_n=1,\quad b_n=0.$$
The structure of the spectrum for such operators is well
known: $\sigma(J) = [-2,2]\bigcup \sigma_d(J)$, where
$\sigma_d(J)$ is at most denumerable set on the complex
plane $\mathbb{C}$ with all accumulation points in
$[-2,2]$. We refer to this portion of spectrum as the {\it
discrete spectrum} of $J$. The goal of our note is to
single out domains on $\mathbb{C}$ free from the discrete
spectrum. In particular, a condition on the matrix entries
which provides the lack of $\sigma_d(J)$ comes in quite
naturally. In the case of selfadjoint operators (\ref{1.3})
($a_n=c_n>0,\ b_n=\overline b_n$) the problem is well
elaborated and goes back to M.S. Birman and J. Schwinger
(see also \cite{G1}, \cite{G2}, \cite{HS} and references
therein, and \cite{L1}, \cite{L2} for complex Jacobi
matrices). The method applied here is totally different and
quite elementary. It is adopted from \cite{M} and based on
certain lower bounds for the Jost function in the unit
disk. The point is made upon the spectrum inclusion regions
rather than bounds for the spectral radius.

{\bf Reccurence relations.} We start out with the three-term
recurrence relation associated with matrix $J$
\begin{equation}
\label{1.6} a_{m-1} y_{m-1} + b_m y_m + c_m y_{m+1} = (z+z^{-1})
y_m,\quad m\in \mathbb{N}, \ \ z\in\overline{\mathbb{D}}:=\{|z|\leq
1\},
\end{equation}
$z\neq 0$ (we put $a_0=c_0 = 1$). It is clear that the initial
data $\{ y_0, y_1\}$ enable one to restore the whole solution
$\{y_m(z)\}_{m\geq 0}$ of (\ref{1.6}), that is, the dimension of
the space of solutions is 2. Sometimes it is beneficial to
deal with a slightly modified relation. If we multiply (\ref{1.6})
by $k(j)=\prod_{i=j}^\infty a_i$ (the product will always be
assumed to converge) and put $x_m = k(m) y_m$, we come to
\begin{equation}\label{1.7}
x_{m-1} + b_m x_m + a_m c_m x_{m+1} = (z + z^{-1}) x_m, \quad m\in
\mathbb{N}.
\end{equation}

{\bf From recurrence relations to discrete integral equations.}
The key role in what follows is played by certain solutions of
(\ref{1.6}), (\ref{1.7}) that have a specific behavior at
infinity. We show that such solutions exist as long as the
coefficients in (\ref{1.6}) tend to their limits fast enough.

Denote by $G$ the Green kernel

\begin{equation}
\label{1.8}
 G(n,m;z)=\left\{ \begin{array}{cc}
\frac{z^{m-n}-z^{n-m}}{z-z^{-1}},&m>n,\\0,&m\leq n,
\end{array}\right.\quad n,m\in\mathbb{Z}_+:=\{0,1,...\},\ z\neq 0.
\end{equation}
It is clear that $$G(n,m,z) =
U_{m-n-1}\left(\frac{z+z^{-1}}{2}\right),$$ where $U_k$ is the
Chebyshev polynomial of the second kind. The recurrence relations
for $G$ are straightforward

\begin{equation}G(n,m+1;z)+G(n,m-1;z)-(z+z^{-1})
G(n,m;z)=\delta(n,m), \label{1.9}
\end{equation}

\begin{equation}G(n-1,m;z)+G(n+1,m;z)-(z+z^{-1})
G(n,m;z)=\delta(n,m), \label{1.10}
\end{equation}
where $\delta(n,m)$  is the Kronecker symbol.

We begin with the following conditional result.

\begin{remark} Suppose that equation $(\ref{1.7})$ has a solution
$v_n$ with asymptotic behavior at infinity
\begin{equation}\label{1.11}
\lim_{n\to\infty}\,v_n(z) z^{-n} = 1
\end{equation}
for some $z\in\mathbb{D}$. Then $v_n$ satisfies a discrete
integral equation
\begin{equation}
\label{1.12} v_n(z)=z^n +\sum_{m=n+1}^\infty J(n,m;z)\,
v_m(z),\quad n\in\mathbb{Z}_+,
\end{equation}
with \begin{equation}\label{1.1}
J(n,m;z)=-b_m G(n,m;z)+\left(1 -
a_{m-1}c_{m-1}\right)G(n,m-1;z).\end{equation}
\end{remark}

{\it Proof.} Let us multiply (\ref{1.9}) by $v_m$,
(\ref{1.7}) by $G(n,m)$ and subtract the latter from the former:
$$G(n,m+1)v_m +G(n,m-1)v_m - G(n,m)v_{m-1}-b_mG(n,m)v_m - a_m c_m
G(n,m)v_{m+1} = \delta(n,m)v_m.$$ Summing up over $m$ from $n$ to
$N$ gives
 $$ v_n = \sum_{m=n}^N\,\left\{-b_m G(n,m)+\left(1 -
a_{m-1}c_{m-1}\right)G(n,m-1)\right\}v_m +$$ $$+
G(n,N+1)v_N - a_Nc_NG(n,N)v_{N+1}. $$ For $|z|<1$ we have
by (\ref{1.8}) and (\ref{1.11}) $$\lim_{N\to
\infty}\,(G(n,N+1)v_N - a_n c_n G(n,N)v_{N+1}) = z^n,$$
which along with $J(n,n)=0$ leads to (\ref{1.12}), as
needed.\hfill $\square$

The converse statement is equally simple.

\begin{remark} Each solution $\{v_n(z)\}_{n\geq0}$,
$z\in\overline{\mathbb{D}}$, of equation $(\ref{1.12})$
satisfies the three-term recurrence relation $(\ref{1.7})$.
\end{remark}

{\it Proof.} Write for $n\geq 1$ $$v_{n-1} + v_{n+1} =
z^{n-1} + z^{n+1} + \sum_{m=n}^\infty\,J(n-1,m)v_m
+\sum_{m=n+2}^\infty\,J(n+1,m)v_m=$$ $$(z+z^{-1})z^n + J(n-1,n)
v_n + J(n-1,n+1)v_{n+1} +\sum_{m=n+2}^\infty\,\{J(n-1,m)
+J(n+1,m)\}v_m.$$
By (\ref{1.8}), (\ref{1.1}) and (\ref{1.10})

$$J(n-1,n;z)=-b_n, \qquad J(n-1,n+1;z)=-(z+z^{-1})b_{n+1}+1-a_nc_n$$
and
$$J(n-1,m;z)+J(n+1,m;z)=(z+z^{-1})J(n,m;z). $$
Hence
 $$v_{n-1}+v_{n+1} + b_nv_n -(1 - a_n c_n)v_{n+1} = (z+z^{-1})
 \left(z^n +\sum_{m=n+1}^\infty\,J(n,m)v_m\right) = (z+z^{-1})v_n,$$
which is exactly (\ref{1.7}).\hfill $\square$

{\bf The Jost solution.} To analyze equation (\ref{1.12})
it seems reasonable to introduce new variables $$\tilde
v_n(z):=v_n z^{-n},\quad \tilde
J(n,m;z):=J(n,m;z)z^{m-n},$$ so that
\begin{equation}\label{1.2}\tilde v_n(z) = 1 +
\sum_{m=n+1}^\infty\,\tilde J(n,m;z)\tilde v_m(z),\quad
n\in\mathbb{Z}_+.\end{equation} Now $\tilde J(n,m;\cdot)$ is a
polynomial and since
 $$|G(n,m,z)z^{m-n}|=\frac{|z^{2(m-n)}-1|}{|z-z^{-1}|}
\leq |z|\min\left\{|m-n|,\frac2{|z^2-1|}\right\},$$ the kernel
$\tilde J$ is bounded by
\begin{equation}\label{1.4} |\tilde
J(n,m;z)|\leq |z| d_m \min\left\{|m-n|,\frac2{|z^2-1|}\right\},\quad
d_m:=|b_m|+|1-a_{m-1}c_{m-1}|, \quad z\in\overline{\mathbb{D}}.
\end{equation}

The main result concerning equation (\ref{1.12}) is the
following

\begin{theorem} $(i)$ Suppose that
\begin{equation}\label{1.13}
\sum_{m=1}^\infty d_m\,<\,\infty.
\end{equation}
Then equation $(\ref{1.12})$ has a unique solution $v_n$
such that $v_n$ is analytic in $\mathbb{D}$, continuous on
$\mathbb{D}_1:=\overline{\mathbb{D}}\setminus\{\pm1\}$ and
\footnote{Following the terminology of selfadjoint case, we
call this solution the {\it Jost solution}. The function
$v_0$ is known as the {\it Jost function}.}
\begin{equation}\label{1.14}
|v_n - z^n|\leq|z|^n\left\{\frac{2|z|}{|z^2 -
1|}\sum_{m=n+1}^\infty\,d_m\right\}\,\exp\left\{\frac{2|z|}{|z^2
- 1|}\sum_{m=n+1}^\infty\,d_m\right\},\ \
z\in\mathbb{D}_1,\ n\in\mathbb{Z}_+.
\end{equation}
$(ii)$ Suppose that \begin{equation}\label{1.15} \sum_{m=1}^\infty
m d_m\,<\,\infty.
\end{equation}
Then $v_n$ is analytic in $\mathbb{D}$, continuous on
$\overline{\mathbb{D}}$ and
\begin{equation}\label{1.16}
|v_n - z^n|\leq|z|^n\left\{\sum_{m=n+1}^\infty\,m
d_m\right\}\,\exp\left\{\sum_{m=n+1}^\infty\,m
d_m\right\},\ \ z\in\mathbb{D},\quad n\in\mathbb{Z}_+.
\end{equation}
\end{theorem}

{\it Proof.} The method of successive approximations does
the job. Write (\ref{1.2}) as \begin{equation} \label{1.17} f_n(z)
= g_n(z) +\sum_{m=n+1}^\infty\,\tilde J(n,m;z)f_m(z)
\end{equation}
with
\begin{equation}\label{1.18}
f_m(z):=\tilde v_m(z) - 1,\quad g_n(z):=\sum_{m=n+1}^\infty
\,\tilde J(n,m;z).
\end{equation}
\noindent (i) Put $\sigma_0(n):=\sum_{m=n+1}^\infty\,d_m$,
$\phi(z):=2|z||z^2-1|^{-1}$ and apply (\ref{1.4}) in the
form
\begin{equation}\label{1.19}
|\tilde J(n,m;z)|\leq\phi(z)\,d_m,\quad
z\in\mathbb{D}_1.\end{equation} Then the series in (\ref{1.18})
converges uniformly on compact subsets of $\mathbb{D}_1$ and so
$g_n$ is analytic in $\mathbb{D}$ and continuous on
$\mathbb{D}_1$. Let us begin with $f_{n,1} = g_n$ and denote
$$f_{n,j+1}(z):=\sum_{m=n+1}^\infty\,\tilde J(n,m;z)f_{m,j}(z).$$
We prove by induction starting with $j=1$ that
\begin{equation}\label{1.20}
\left |f_{n,j}(z)\right
|\leq\,\frac{(\phi(z)\sigma_0(n))^j}{(j-1)!}.
\end{equation}
It is obvious for $j=1$ by (\ref{1.19}). Next, let (\ref{1.20}) be
true. Then
$$|f_{n,j+1}(z)|\leq\phi(z)\,\sum_{m=n+1}^\infty\,d_m|f_{m,j}(z)|\leq
\frac{(\phi(z))^{j+1}}{(j-1)!}\sum_{m=n+1}^\infty\,d_m\sigma_0^j(m).$$
An elementary inequality $(a+b)^{j+1} - a^{j+1}\geq(j+1)ba^j$
gives
$$\sum_{m=n+1}^\infty\,d_m\sigma_0^j(m)\leq\frac{1}{j}
\sum_{m=n+1}^\infty\{\sigma_0^{j+1}(m-1)
-\sigma_0^{j+1}(m)\}=\frac{\sigma_0^{j+1}(n)}{j},$$ which
proves (\ref{1.20}) for $f_{n,j+1}$. Thereby the series
$$f_n(z)=\sum_{j=1}^\infty\,f_{n,j}(z)$$
converges uniformly on compact subsets of $\mathbb{D}_1$ and
solves (\ref{1.17}), being analytic in $\mathbb{D}$ and continuous
on $\mathbb{D}_1$. The estimate (\ref{1.14}) follows from (\ref{1.20})
and $\tilde v_n=v_n z^{-n}.$

Suppose that there are two solutions $f_n$ and $\tilde f_n$
of (\ref{1.17}). Take the difference and apply
(\ref{1.19}):
\begin{equation}
h_n\leq \sum_{m=n+1}^\infty \phi(z) h_md_m=q_n, \qquad
h_n:=|f_n(z) -\tilde f_n(z)|. \label{1.50}
\end{equation}
Clearly, $q_n\to 0$ as $n\to\infty$ and $q_k=0$ for some
$k$ implies by (\ref{1.50}) $h_n\equiv 0$. If $q_n>0$, then
\begin{equation}
\frac{q_{n-1}-q_n}{q_n}=\frac{h_n\phi(z)d_n}{q_n}\leq
\phi(z)d_n, \quad
q_k\leq\prod_{j=k+1}^M\left(1+\phi(z)d_j\right)\,q_M
\end{equation}
which leads to $q_k=0$ and again $h_n\equiv 0$. So the
uniqueness is proved.

\noindent
(ii) The same sort of reasoning is applicable
with $$|\tilde J(n,m;z)|\leq|z||m-n|d_m\leq m d_m$$ and
$$|f_{n,j}(z)|\leq\frac{\sigma_1^j(n)}{(j-1)!},\quad
\sigma_1(n):=\sum_{m=n+1}^\infty m d_m$$ instead of (\ref{1.19})
and (\ref{1.20}), respectively.\hfill $\square$
\begin{corollary}
\label{cor1} Let $t$ be the root of the equation \begin{equation}
te^t=1,\quad t\approx 0.567.\label{1.21}
\end{equation}
Under assumption $(\ref{1.13})$ the Jost function $v_0$ does not
vanish in the domain
\begin{equation}
\Omega:=\{z\in\mathbb{D};\quad |z-z^{-1}|>2t^{-1}\sum_{m=1}^\infty
d_m\}.\label{1.22}\end{equation} Under assumption
$(\ref{1.15})$
$v_0$ does not vanish in $\mathbb{D}$ as long as
\begin{equation}
\sum_{m=1}^\infty m d_m <t.\label{1.23}\end{equation}
\end{corollary}

{\bf Eigenvalues of the Jacobi matrix and zeros of the Jost
function}. Going back to the matrix $J$, let
$\lambda\in\sigma_d(J)$ with an eigenvector
$h=\{h_n\}_{n\geq 1}$. Clearly, $h_1\neq 0$. As a result
there is only one linear independent eigenvector
corresponding to any eigenvalue $\lambda$ (an appropriate
linear combination of two would give zero). Denote by
$$Z=Z(J) = \{z\in \mathbb{D}:\, v_0(z)=0\}$$ the zero set of
the Jost function $v_0$ in $\mathbb{D}$. As $\{v_n\}_{n\geq
0}\in\ell^2$, for each $z_0\in Z$ the vector
$\{v_n(z_0)\}_{n\geq 1}$ is the eigenvector of $J$ with the
eigenvalue $\lambda_0=z_0 +z_0^{-1}.$ Conversely, let
$\{g_n\}$ and $\{h_n\}$ be two solutions of (\ref{1.6}).
Define their "wronskian" by $$W_n(g,h)=g_nh_{n+1} -
g_{n+1}h_n,\quad n\geq 0.$$ From (\ref{1.6}) it follows
that $a_{n-1}W_{n-1} = c_n W_n$, and iterating gives
\begin{equation}
\label{1.24} W_n=\,\frac{a_{n-1}...a_0}{c_nc_{n-1}...c_1}\,
W_0.\end{equation} Suppose that both solutions are from $\ell^2$.
Then $W_n$ goes to zero as $n\to\infty$ and since both products in
the RHS of (\ref{1.24}) converge, $W_0=0$. The latter means that
$g$ and $h$ are linearly dependent, and so $u_0(z_0) = 0.$ We end
up with
\begin{remark}
 $$\sigma_d(J) = \{z+z^{-1}:\ \ z\in Z(J)\}.$$
\end{remark}
Our main result concerning the discrete spectrum of the Jacobi matrices
(\ref{1.3}) can be displayed as follows.
\begin{theorem}
Under assumption $(\ref{1.13})$ the domain $$G(J)=\{z+z^{-1}:\ \
z\in\Omega\}$$ with $\Omega$ $(\ref{1.22})$ is free from the
discrete spectrum $\sigma_d(J)$. The matrix $J$ has no discrete
spectrum at all as soon as $(\ref{1.23})$ holds, where $t$ is the
solution of $(\ref{1.21})$.
\end{theorem}
{\bf Remark}. Suppose that
$$c=\frac{2}{t}\,\sum_{m=1}^\infty d_m<2. $$ Then
$\sigma_d(J)$ is contained in the union of two symmetric
rectangles
$$\sigma_d(J)\subset\left\{w:\,\sqrt{4-c^2}<|\mbox{Re}\,w|<\sqrt{4+c^2},\
\ |\mbox{Im}\,w|<\frac{c^2}{4}\right\}. $$

\end{document}